# Statistical analysis of redundant systems with "warm" stand-by units

V. Bagdonavičius[a]; I. Masiulaitytė[a] and M. Nikulin[b]

[a] *University of Vilnius, Lithuania*; [b] *Victor Segalen University, Bordeaux, France*

Mathematical formulation of fluent switching from "warm" to "hot" conditions of standby units is given using the well known Sedyakin's and accelerated failure time (AFT) models. Non-parametric estimators of cumulative distribution function and mean failure time of a redundant system with several stand-by units are proposed. Goodness-of-fit tests for two given models are given.

**Keywords:** AFT model; Chi-squared goodness-of-fit; Non-parametric estimation; Reliability; Redundant system; Sedyakin model, Stand-by

**AMS Subject Classification**: 60K10; 62G05; 62G10; 62G20; 62G30

## 1. Introduction

To warrant high reliability of key components of systems, stand-by units are used. If such a component fails and at least one not failed stand-by unit remains, it operates instead of the key component. If this one fails, the second stand-by unit is used and so on until the last stand-by unit fails. We suppose that switching is instantaneous and there are no repairs.

If the stand-by units are functioning in the same "hot" conditions as the main unit then usually after switching the reliability of the stand-by units does not change. But "hot" redundancy has disadvantages because any one of the stand-by units fails earlier than the main one with the probability 0.5.

If the stand-by units are not operating until the failure of the main unit ("cold" reserving), it is possible that during and after switching the failure rate increases because the stand-by unit is not "warmed" enough [1]. So "warm" reserving is sometimes used [2]: stand by units function under lower stress than the main one. In such a case the probability of the failure of the stand-by unit is smaller than that of the main unit and it is also possible that switching is fluent, i.e. switching from "warm" to "hot" conditions does not do any damage to units. What does it mean mathematically?

Let us consider a system of $m$ units: one main unit and $m-1$ stand-by units. We shall use notation $S(1, m-1)$ for such systems.

Denote by $T_1$, $F_1$ and $f_1$ the failure time, the c.d.f. and the probability density function of the main unit. The failure times of the stand-by units are denoted by $T_2, \ldots, T_m$. In "hot" conditions their distribution functions are also $F_1$. In "warm" conditions the c.d.f. of $T_i$ is $F_2$ and the p.d.f is $f_2$, $i = 2, \ldots, m$. If a stand-by unit is switched from "warm" to "hot" conditions, its c.d.f. is different from $F_1$ and $F_2$.

In 1966 N.M.Sedyakin [3] formulated the "physical principle in reliability". The idea is the following. For two groups of units from identical populations functioning under different constant stresses $x_1$ and $x_2$, two times $t_1$ and $t_2$ are equivalent in the sense that the probabilities of survival until these moments are equal:

$$1 - F_1(t_1) = \mathbf{P}(T_1 > t_1|x_1) = \mathbf{P}(T_2 > t_2|x_2) = 1 - F_2(t_2). \tag{1}$$

If the first group of units is tested under the constant stress $x_1$, the second group is tested under the constant stress $x_2$ up to time $t_2$ and under the stress $x_1$ for $t \geq t_2$, i.e. under the step-stress

$$x(u) = \begin{cases} x_2, & 0 \leq u < t_2, \\ x_1, & u \geq t_2, \end{cases}$$

then Sedyakin's hypothesis is the following: for all $s > 0$

$$\lambda(t_2 + s|x(\cdot)) = \lambda(t_1 + s|x_1); \tag{2}$$

here $\lambda(t|x)$ denotes the hazard function under the stress $x$ and $t_1$ satisfies (1). The equality implies

$$\mathbf{P}(T_2 > t_2 + s|T_2 > t_2, x(\cdot)) = \mathbf{P}(T_1 > t_1 + s|T_1 > t_1, x_1),$$

which means that the probability not to fail time $s$ under the same stress after two equivalent times $t_1$ and $t_2$ are the same for two groups of units. So switching from the stress $x_1$ to the stress $x_2$ does not do any damage to units. If the stresses $x_2$ and $x_1$ mean "warm" and "hot" conditions, respectively, and the switch on time $t_2$ would be non-random then we could use formula (2) for computing the conditional distribution of the switched on stand-by unit (see [6],[7]). In the case of redundant systems the formula (2) should be modified because the switch on times are random - stand-by units are switched on after the failure of the unit operating in "hot" conditions.

The failure time of the system $S(1, m-1)$ is $T^{(m)} = T_1 \vee T_2 \vee \ldots \vee T_m$. Using the principle of Sedyakin, we consider the following model for the c.d.f. of the random variable $T$. As $T^{(m)} = (T_1 \vee T_2 \vee \ldots \vee T_{m-1}) \vee T_m$, we can consider this system as a system $S(1,1)$ with one main element (which itself is a system $S(1, m-2)$) and one stand-by element.

Denote by $K_j$ and $k_j$ the c.d.f. and the p.d.f. of $T^{(j)}$, respectively, $(j = 2, \ldots, m)$, $K_1 = F_1$, $k_1 = f_1$. The c.d.f $K_j$ can be written in terms of the c.d.f $K_{j-1}$ and $F_1$:

$$K_j(t) = \mathbf{P}(T^{(j)} \leq t) = \mathbf{P}(T^{(j-1)} \leq t, T_j \leq t) = \int_0^t \mathbf{P}(T_j \leq t|T^{(j-1)} = y)dK_{j-1}(y). \tag{3}$$

We use the Sedyakin's principle modelling the conditional distribution $\mathbf{P}(T_j \leq t|T^{(j-1)} = y)$. If $T^{(j-1)} = y$ and the stand-by unit has not failed until this moment then it is switched on.

Hypothesis $H_0$:

$$f_{T_j|T^{(j-1)}=y}(t) = \begin{cases} f_2(t) & \text{if } t \leq y, \\ f_1(t + g(y) - y) & \text{if } t > y; \end{cases} \tag{4}$$

"equivalent to $y$ moment" $g(y)$ is found from the equation $F_1(g(y)) = F_2(y)$, so

$$g(y) = F_1^{-1}(F_2(y)).$$

The formula (4) implies

$$K_j(t) = \int_0^t F_1(t + g(y) - y)dK_{j-1}(y). \tag{5}$$

So the distribution function of the system with $m-1$ stand-by units is defined recurrently using formula (5) ($j = 2, \ldots, m$).

In particular, if we suppose that the distribution of units functioning in "warm" and "hot" conditions differ only in scale, i.e.

$$F_2(t) = F_1(rt), \tag{6}$$

for all $t \geq 0$ and some $r > 0$, then $g(y) = ry$.

Hypothesis $H_0^*$:

$$f_{T_j|T^{(j-1)}=y}(t) = \begin{cases} f_2(t) & \text{if } t \leq y, \\ f_1(t + ry - y) & \text{if } t > y., \end{cases} \tag{7}$$

Conditionally (given $T^{(j-1)} = y$) this hypothesis corresponds to the accelerated failure time (AFT) model [4], [5], [8]. Under model (7) the distribution function of the system is obtained using recurrent formulas

$$K_j(t) = \int_0^t F_1(t + ry - y)dK_{j-1}(y). \tag{8}$$

## 2. Nonparametric estimation of the reliability of the redundant system

Suppose that the hypothesis $H_0^*$ is true and the following data are available :

a) complete ordered sample $T_{11}, \ldots, T_{1n_1}$ of the failure times of units tested in "hot" conditions;

b) the time to obtain complete data in "warm" conditions may be long, so we suppose that $n_2$ units are tested up to time $t_1$ in "warm" conditions and the ordered first failure times $T_{21}, \ldots, T_{2m_2}$ are obtained.

Set

$$N_1(t) = \sum_{i=1}^{n_1} \mathbf{1}_{\{T_{1i} \leq t\}}, \quad N_2(t) = \sum_{i=1}^{n_2} \mathbf{1}_{\{T_{2i} \leq t, t \leq t_1\}},$$

$$Y_1(t) = \sum_{i=1}^{n_1} \mathbf{1}_{\{T_{1i} \geq t\}}, \quad Y_2(t) = \sum_{i=1}^{n_2} \mathbf{1}_{\{T_{2i} \geq t, t \leq t_1\}}.$$

Note that the random variables $T_{1i}/r$ and $T_{2i}$ can be interpreted as order statistics from samples of size $n_1$ and $n_2$, respectively, from the population having the c.d.f

$F_2$. So if we denote

$$\tilde{N}_1(t) = \sum_{i=1}^{n_1} \mathbf{1}_{\{T_{1i}/r \leq t\}} = N_1(rt), \quad \tilde{N}_2(t) = N_2(t),$$

$$\tilde{Y}_1(t) = \sum_{i=1}^{n_1} \mathbf{1}_{\{T_{1i}/r \geq t\}} = Y_1(rt), \quad \tilde{Y}_2(t) = Y_2(t),$$

then Nelson-Aaalen type estimator (still depending on $r$) of the cumulative hazard function $\Lambda_2 = -\ln S_2$ can be considered:

$$\tilde{\Lambda}_2(t, r) = \int_0^t \frac{d\tilde{N}_1(u) + d\tilde{N}_2(u)}{\tilde{Y}_1(u) + \tilde{Y}_2(u)} = \int_0^t \frac{dN_1(ru) + dN_2(u)}{Y_1(ru) + Y_2(u)}. \tag{9}$$

Taking into consideration that the difference

$$M_2(t) = N_2(t) - \int_0^t Y_2(u) d\Lambda_2(u)$$

is a martingale on $[0, t_1]$ with respect to the filtration generated by the data, and $\mathbf{E}M_2(t_1) = 0$, the parameter $r$ can be estimated using the estimating function

$$U(r) = N_2(t_1) - \int_0^{t_1} Y_2(u) d\tilde{\Lambda}_2(u, r) =$$

$$N_2(t_1) - \int_0^{rt_1} \frac{Y_2(v/r) dN_1(v)}{Y_1(v) + Y_2(v/r)} - \int_0^{t_1} \frac{Y_2(u) dN_2(u)}{Y_1(ru) + Y_2(u)}.$$

$U(r)$ is a non-increasing step function,

$$U(0+) = N_2(t_1) - \int_0^{t_1} \frac{Y_2(u) dN_2(u)}{n_1 + Y_2(u)} > 0, \quad U(+\infty) = -\int_0^{\infty} \frac{n_2 dN_1(v)}{Y_1(v) + n_2} < 0,$$

so the parameter $r$ is estimated by the statistic

$$\hat{r} = U^{-1}(0) = sup\{r : U(r) > 0\}.$$

The estimator of the cumulative hazards $\Lambda_1$ and $\Lambda_2$ are

$$\hat{\Lambda}_1(t) = \tilde{\Lambda}_2(t/\hat{r}, \hat{r}) = \int_0^t \frac{dN_1(u)}{Y_1(u) + Y_2(u/\hat{r})} + \int_0^{t/\hat{r}} \frac{dN_2(u)}{Y_1(\hat{r}u) + Y_2(u)} =$$

$$\sum_{T_{1i} \leq t} \frac{1}{Y_1(T_{1i}) + Y_2(T_{1i}/\hat{r})} + \sum_{T_{2i} \leq t/\hat{r}} \frac{1}{Y_1(\hat{r}T_{2i}) + Y_2(T_{2i})},$$

$$\hat{\Lambda}_2(t) = \hat{\Lambda}_1(\hat{r}t).$$

The estimators of the c.d.f. $S_i = 1 - F_i$ are the product integrals of the estimators $\hat{\Lambda}_1(t)$ and $\hat{\Lambda}_2(t)$, so

$$\hat{F}_1(t) = 1 - \pi_{0 \leq s \leq t}(1 - d\hat{A}_1(s)) =$$

$$1 - \prod_{T_{1i} \leq t}\left(1 - \frac{1}{Y_1(T_{1i}) + Y_2(T_{1i}/\hat{r})}\right) \prod_{T_{2i} \leq t/\hat{r}}\left(1 - \frac{1}{Y_1(\hat{r}T_{2i}) + Y_2(T_{2i})}\right),$$

$$\hat{F}_2(t) = \hat{F}_1(\hat{r}t).$$

Mixing all moments $T_{1i}$ and $\hat{r}T_{2j}$ and ordering them, we obtain the sequence of random variables $T_1 \leq \ldots \leq T_{n_1+m_2}$. The estimators $\hat{F}_1(t)$ and $\hat{F}_2(t)$ can be written:

$$\hat{F}_1(t) = 1 - \prod_{T_i \leq t}\left(1 - \frac{1}{Y_1(T_i) + Y_2(T_i/\hat{r})}\right), \quad \hat{F}_2(t) = \hat{F}_1(\hat{r}t).$$

The c.d.f. $K_m$ of the redundant system is estimated using the following recurrent equations $(j = 2, \ldots, m)$:

$$\hat{K}_j(t) = \hat{F}_2(t)\hat{K}_{j-1}(t) + \int_{\hat{r}t}^{t} \hat{K}_{j-1}\left(\frac{t-z}{1-\hat{r}}\right) d\hat{F}_1(z) =$$

$$\hat{F}_2(t)\hat{K}_{j-1}(t) + \sum_{\hat{r}t < T_i \leq t} \hat{K}_{j-1}\left(\frac{t - T_i}{1 - \hat{r}}\right) \frac{\hat{F}_1(T_{i-1})}{Y_1(T_i) + Y_2(T_i/\hat{r})}. \tag{10}$$

The estimator of the mean failure time $\mu$ of the system is

$$\hat{\mu} = \int_0^\infty t d\hat{K}_m(t) = \sum_{i=1}^{n_1+m_2} T_i [\hat{K}_m(T_i) - \hat{K}_m(T_{i-1})].$$

Failure times from exponential distribution were simulated:

$$T_{1j} \sim E(\lambda_1), \quad T_{2j} \sim E(\lambda_2), \quad \lambda_2 = r\lambda_1.$$

The graphs of the estimators of the c.d.f. $F_1$ and $K_m$ ($m = 2, 3, 4$), in the case of complete samples and different sample sizes are presented in Figures 1 and 2. Increasing the number of stand-by units increases the reliability of the redundant system.

## 3. Goodness-of-fit

The given estimators $\hat{K}_m$ of the c.d.f. of the redundant system $S(1, m-1)$ can be used if the hypothesis $H_0$ is true. If switching from "warm" to "hot" conditions does not damage units in the system $S(1,1)$ then it is natural that this is true for the system $S(1, m-1)$, $m > 2$. So we do not need a test with $S(1, m-1)$ systems and it is sufficient to give tests for the hypothesis $H_0$ and $H_0^*$ when only one stand-by unit is used.

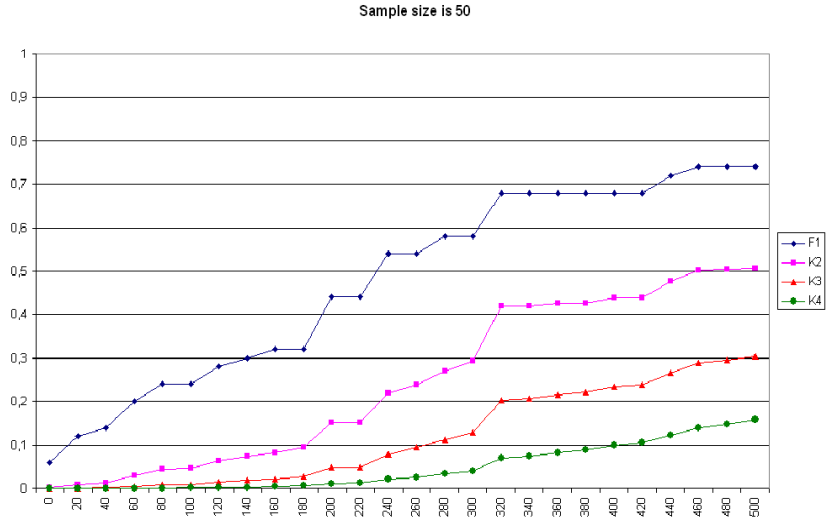

Figure 1. Estimators of the c.d.f. $F_1$ and $K_m$ $(m = 2, 3, 4)$, $n_1 = n_2 = 50$, $t_1 = \infty$

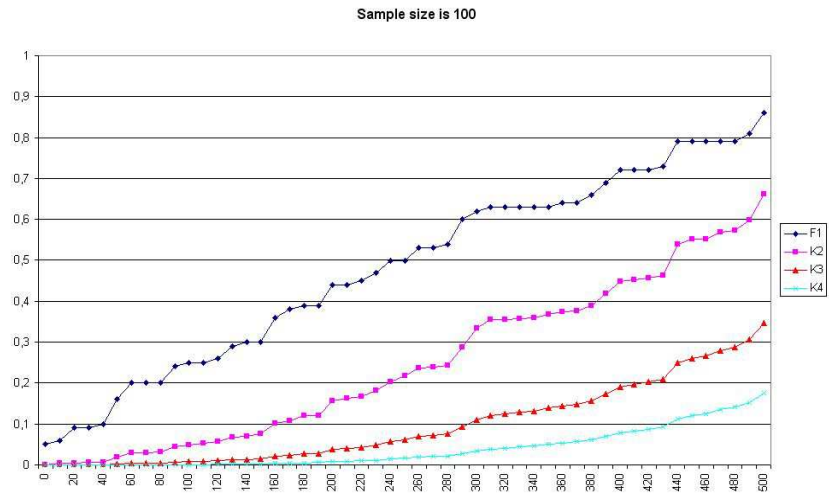

Figure 2. Estimators of the c.d.f. $F_1$ and $K_m$ $(m = 2, 3, 4)$, $n_1 = n_2 = 100$, $t_1 = \infty$

Suppose that the following data are available :

a) the failure times $T_{11}, \ldots, T_{1n_1}$ of $n_1$ units tested in "hot" conditions;

b) the failure times $T_{21}, \ldots, T_{2n_2}$ of $n_2$ units tested in "warm" conditions;

c) the failure times $T_1, \ldots, T_n$ of $n$ redundant systems (with "warm" stand-by units).

The tests are based on the difference of two estimators of the c.d.f. $F$. The first estimator is the empirical distribution function

$$\hat{F}^{(1)}(t) = \frac{1}{n} \sum_{i=1}^{n} \mathbf{1}_{\{T_i \leq t\}}. \tag{11}$$

The second is based on formula (5), i.e.

$$\hat{F}^{(2)}(t) = \int_0^t \hat{F}_1(t + \hat{g}(y) - y) d\hat{F}_1(y),$$

where (hypothesis $H_0$)

$$\hat{g}(y) = \hat{F}_1^{-1}(\hat{F}_2(y)), \quad \hat{F}_j(t) = \frac{1}{n_j} \sum_{i=1}^{n_j} \mathbf{1}_{\{T_{ji} \leq t\}}, \quad \hat{F}_1^{-1}(y) = \inf\{s : \hat{F}_1(s) \geq y\}, \tag{12}$$

or (hypothesis $H_0^*$)

$$\hat{g}(y) = \hat{r}y, \quad \hat{r} = \frac{\hat{\mu}_1}{\hat{\mu}_2}, \quad \hat{\mu}_j = \frac{1}{n_j} \sum_{i=1}^{n_j} T_{ji}. \tag{13}$$

The test is based on the statistic

$$X = \sqrt{n} \int_0^\infty [\hat{F}^{(1)}(t) - \hat{F}^{(2)}(t)] dt. \tag{14}$$

It is a natural generalization of Student's t-test for comparing the means of two populations. Indeed, the mean failure time of the system with c.d.f. $F$ is

$$\mu = \int_0^\infty [1 - F(s)] ds,$$

so the statistic (12) is the normed difference of two estimators (the second is not the empirical mean) of the mean $\mu$. Student's t-test is based on the difference of empirical means of two populations.

It will be shown that in the case of both hypothesis $H_0$ and $H_0^*$ the limit distribution (as $n_i/n \to l_i \in (0,1)$, $n \to \infty$) of the statistic $X$ is normal with zero mean and finite variance $\sigma^2$, (see Theorems 1 and 2).

The test statistic is

$$Y_n^2 = \left(\frac{X}{\hat{\sigma}}\right)^2$$

where $\hat{\sigma}$ is a consistent estimator of $\sigma$. The distribution of the statistic $Y_n$, is approximated by the standard normal distribution and the hypothesis $H_0$ (or $H_0^*$) is rejected with approximative significance value $\alpha, 0 < \alpha < 0.5$, if $Y_n^2 > \chi_{1-\alpha}^2(1)$, where $\chi_{1-\alpha}^2(1)$ is the $(1-\alpha)$-quantile of the chi-squared distribution with one degree of freedom, see, for example, Greenwood and Nikulin (1996).

Let us find the asymptotic distribution of the statistic (14).

**Theorem 3.1:** *Suppose that $n_i/n \to l_i \in (0,1)$, $n \to \infty$ and the densities $f_i(x)$, $i = 1, 2$ are continuous and positive on $(0, \infty)$. Then under $H_0^*$ the statistic (14) converges in distribution to the normal law $N(0, \sigma^2)$, where*

$$\sigma^2 = \mathbf{Var}(T_i) + \frac{1}{l_1} \mathbf{Var}(H(T_{1i})) + \frac{c^2 r^2}{l_2} \mathbf{Var}(T_{2i}), \tag{15}$$

*where*

$$H(x) = x[c + r - 1 - F_1(x/r) - rF_2(x)] + r\mathbf{E}(\mathbf{1}_{\{T_{1i} \leq x/r\}} T_{1i}) + r\mathbf{E}(\mathbf{1}_{\{T_{2i} \leq x\}} T_{2i}),$$

$$c = \frac{1}{\mu_2} \int_0^\infty y[1 - F_2(y)]dF_1(y).$$

**Proof:** The limit distribution of the empirical distribution functions is well known:

$$\sqrt{n}(\hat{F}_i - F_i) \xrightarrow{\mathcal{D}} U_i, \quad \sqrt{n}(\hat{F}^{(1)} - F) \xrightarrow{\mathcal{D}} U \qquad (16)$$

on $D[0, \infty)$, where $\xrightarrow{\mathcal{D}}$ means weak convergence, $U_1, U_2$ and $U$ are independent Gaussian martingales with $U_i(0) = U(0) = 0$ and the covariances

$$\mathbf{cov}(U_i(s_1), U_i(s_2)) = \frac{1}{l_i} F_i(s_1 \wedge s_2)[1 - F_i(s_1 \vee s_2)],$$

$$\mathbf{cov}(U(s_1), U(s_2)) = F(s_1 \wedge s_2)[1 - F(s_1 \vee s_2)].$$

Under hypothesis $H_0^*$ the difference of the two estimators of the distribution function $F$ can be written as follows:

$$\hat{F}^{(1)}(t) - \hat{F}^{(2)}(t) = \hat{F}^{(1)}(t) - F(t) - \int_0^t \hat{F}_1(t + \hat{g}(y) - y)d\hat{F}_1(y) + \int_0^t F_1(t + g(y) - y) \times$$

$$dF_1(y) = \hat{F}^{(1)}(t) - F(t) - \int_0^t [F_1(t + \hat{g}(y) - y) - F_1(t + g(y) - y)]dF_1(y) -$$

$$\int_0^t [(\hat{F}_1(t + \hat{g}(y) - y) - \hat{F}_1(t + g(y) - y)) - (F_1(t + \hat{g}(y) - y) - F_1(t + g(y) - y))]dF_1(y) -$$

$$\int_0^t [\hat{F}_1(t + \hat{g}(y) - y) - \hat{F}_1(t + g(y) - y)][d\hat{F}_1(y) - dF_1(y)] -$$

$$\int_0^t [\hat{F}_1(t + g(y) - y) - F_1(t + g(y) - y)]dF_1(y) -$$

$$\int_0^t [\hat{F}_1(t + g(y) - y) - F_1(t + g(y) - y)][d\hat{F}_1(y) - dF_1(y)] -$$

$$\int_0^t F_1(t + g(y) - y)[d\hat{F}_1(y) - dF_1(y)].$$

The statistic (14) can be written

$$X = \int_0^\infty \sqrt{n}[\hat{F}^{(1)}(t) - F(t)]dt-$$

$$\int_0^\infty dt \int_0^t \sqrt{n}[F_1(t+\hat{g}(y)-y) - F_1(t+g(y)-y)]\,dF_1(y)-$$

$$\int_0^\infty dt \int_0^t \sqrt{n}[\hat{F}_1(t+g(y)-y) - F_1(t+g(y)-y)]\,dF_1(y)-$$

$$\int_0^\infty dt \int_0^t F_1(t+g(y)-y)d\{\sqrt{n}[\hat{F}_1(y) - F_1(y)]\} + o_P(1). \qquad (17)$$

Set $\sigma_j^2 = \mathbf{Var}(T_{ji})$, $j = 1, 2$. The convergence

$$\sqrt{n}(\hat{\mu}_j - \mu_j) \xrightarrow{D} Y_j = -\int_0^\infty U_j(y)dy \sim N(0, \sigma_j^2/l_i)$$

implies

$$\sqrt{n}(\hat{r} - r) \xrightarrow{D} Y = \frac{1}{\mu_2}(Y_1 - rY_2) \sim N(0, \frac{\sigma_1^2}{\mu_2^2}(\frac{1}{l_1} + \frac{1}{l_2})). \qquad (18)$$

Formulas (16)-(18) imply

$$\int_0^\infty \sqrt{n}[\hat{F}^{(1)}(t) - F(t)]dt \xrightarrow{D} \int_0^\infty U(t)dt,$$

$$\int_0^\infty dt \int_0^t \sqrt{n}[F_1(t+\hat{g}(y)-y) - F_1(t+g(y)-y)]\,dF_1(y) \xrightarrow{D}$$

$$cY_1 - rcY_2 = -c\int_0^\infty U_1(y)dy + rc\int_0^\infty U_2(y)dy,$$

$$\int_0^\infty dt \int_0^t \sqrt{n}[\hat{F}_1(t+g(y)-y) - F_1(t+g(y)-y)]\,dF_1(y) \xrightarrow{D}$$

$$\int_0^\infty dt \int_0^t U_1(t+g(y)-y)dF_1(y) = \int_0^\infty dF_1(y) \int_{g(y)}^\infty U_1(u)du =$$

$$\int_0^\infty U_1(u)F_1(g^{-1}(u))du,$$

$$\int_0^\infty dt \int_0^t F_1(t+g(y)-y)d\{\sqrt{n}[\hat{F}_1(y)-F_1(y)]\} \xrightarrow{D}$$

$$\int_0^\infty dt \int_0^t F_1(t+g(y)-y)dU_1(y) = \int_0^\infty F_2(t)U_1(t)dt-$$

$$\int_0^\infty U_1(y)[1-F_2(y)]d(g(y)-y) = \int_0^\infty U_1(y)[rF_2(y)-r+1]dy.$$

We obtained

$$X \xrightarrow{D} V_1 + V_2 + V_3,$$

where

$$V_1 = \int_0^\infty U(y)dy, \quad V_2 = \int_0^\infty h(y)U_1(y)dy, \quad h(y) = c+r-1-F_1(y/r)-rF_2(y),$$

$$V_3 = -rc\int_0^\infty U_2(y)dy.$$

The variances of the random variables $V_i$ are:

$$\mathbf{Var}(V_1) = \mathbf{Var}(T_i), \quad \mathbf{Var}(V_3) = \frac{c^2 r^2}{l_2}\mathbf{Var}(T_{2i})$$

$$\mathbf{Var}(V_2) = \frac{2}{l_1}\int_0^\infty [1-F_1(y)]h(y)dy \int_0^y F_1(z)h(z)dz = \frac{1}{l_1}\mathbf{Var}(H(T_{1i})),$$

where

$$H(x) = \int_0^x h(y)dy = x[c+r-1-F_1(x/r)-rF_2(x)] + \int_0^x ydF_1(y/r)+$$

$$r\int_0^x ydF_2(y) = x[c+r-1-F_1(x/r)-rF_2(x)]+$$

$$r\mathbf{E}(\mathbf{1}_{\{T_{1i}\leq x/r\}}T_{1i}) + r\mathbf{E}(\mathbf{1}_{\{T_{2i}\leq x\}}T_{2i}).$$

$\square$

A consistent estimator of the variance $\sigma^2$ is

$$\hat{\sigma}^2 = \frac{1}{n}\sum_{i=1}^n (T_i - \hat{\mu})^2 + \frac{n}{n_1^2}\sum_{i=1}^{n_1}[\hat{H}(T_{1i}) - \hat{\bar{H}}]^2 + \frac{\hat{c}^2 \hat{r}^2 n}{n_2^2}\sum_{i=1}^{n_2}(T_{2i} - \hat{\mu}_2)^2,$$

where

$$\hat{\mu} = \frac{1}{n}\sum_{i=1}^{n} T_i, \quad \hat{c} = \frac{1}{\hat{\mu}_2}\int_0^\infty y[1-\hat{F}_2(y)]d\hat{F}_1(y) = \frac{1}{\hat{\mu}_2 n_1}\sum_{i=1}^{n_1} T_{1i}[1-\hat{F}_2(T_{1i})],$$

$$\hat{H}(x) = x[\hat{c}+\hat{r}-1-\hat{F}_1(x/\hat{r})-\hat{r}\hat{F}_2(x)] + \frac{\hat{r}}{n_1}\sum_{i=1}^{n_1}\mathbf{1}_{\{T_{1i}\leq x/\hat{r}\}}T_{1i} + \frac{\hat{r}}{n_2}\sum_{i=1}^{n_2}\mathbf{1}_{\{T_{2i}\leq x\}}T_{2i},$$

$$\hat{\bar{H}} = \frac{1}{n_1}\sum_{i=1}^{n_1}\hat{H}(T_{1i}).$$

**Theorem 3.2:** *Suppose that $n_i/n \to l_i \in (0,1)$, $n \to \infty$ and the densities $f_i(x)$, $i=1,2$ are continuous and positive on $(0,\infty)$. Then under $H_0$ the statistic (14) converges in distribution to the normal law $N(0,\sigma^2)$, where*

$$\sigma^2 = \mathbf{Var}(T_i) + \frac{1}{l_1}\mathbf{Var}(H(T_{1i})) + \frac{1}{l_2}\mathbf{Var}(Q(T_{2i}))$$

where

$$H(x) = Q(x) - xF_1(g^{-1}(x)) + g(x)[1-F_2(x)] + \mathbf{E}(\mathbf{1}_{\{g(T_{1i})\leq x\}}g(T_{1i})) + \mathbf{E}(\mathbf{1}_{\{T_{2i}\leq x\}}g(T_{2i})) - x,$$

$$Q(x) = \mathbf{E}\{\mathbf{1}_{\{T_{1i}\leq x\}}[1-F_2(T_{1i})]/f_1(g(T_{1i}))\}.$$

**Proof:** Similarly as in Theorem 1 we obtain

$$\int_0^\infty \sqrt{n}[\hat{F}^{(1)}(t) - F(t)]dt \xrightarrow{\mathcal{D}} \int_0^\infty U(t)dt,$$

$$\int_0^\infty dt \int_0^t \sqrt{n}[F_1(t+\hat{g}(y)-y) - F_1(t+g(y)-y)]\,dF_1(y) \xrightarrow{\mathcal{D}}$$

$$-\int_0^\infty \frac{U_1(g(y)) - U_2(y)}{f_1(g(y))} f_1(y)[1-F_2(y)]dy,$$

$$\int_0^\infty dt \int_0^t \sqrt{n}[\hat{F}_1(t+g(y)-y) - F_1(t+g(y)-y)]\,dF_1(y) \xrightarrow{\mathcal{D}}$$

$$\int_0^\infty U_1(u)F_1(g^{-1}(u))du,$$

$$\int_0^\infty dt \int_0^t F_1(t+g(y)-y)d\{\sqrt{n}[\hat{F}_1(y) - F_1(y)]\} \xrightarrow{\mathcal{D}}$$

$$\int_0^\infty dt \int_0^t F_1(t+g(y)-y)dU_1(y) = \int_0^\infty F_2(t)U_1(t)dt-$$

$$\int_0^\infty U_1(y)[1-F_2(y)]d(g(y)-y) = \int_0^\infty U_1(y)\{F_2(y)-(g'(y)-1)[1-F_2(y)]\}dy.$$

We obtained

$$X \xrightarrow{D} V_1 + V_2 + V_3,$$

where

$$V_1 = \int_0^\infty U(y)dy, \quad V_2 = \int_0^\infty h(y)U_1(y)dy,$$

$$h(y) = \frac{f_1(y)}{f_1(g(y))}[1-F_2(y)] - F_1(g^{-1}(y)) - F_2(y) + (g'(y)-1)[1-F_2(y)].$$

$$V_3 = -\int_0^\infty \frac{U_2(y)}{f_1(g(y))}[1-F_2(y)]dF_1(y).$$

The variances of the random variables $V_i$ are:

$$\mathbf{Var}(V_1) = \mathbf{Var}(T_i),$$

$$\mathbf{Var}(V_3) = \frac{2}{l_2}\int_0^\infty \frac{[1-F_2(y)]^2 dF_1(y)}{f_1(g(y))} \int_0^y \frac{F_2(z)[1-F_2(z)]dF_1(z)}{f_1(g(z))} = \frac{1}{l_2}\mathbf{Var}(Q(T_{2i})),$$

$$\mathbf{Var}(V_2) = \frac{2}{l_1}\int_0^\infty [1-F_1(y)]h(y)dy \int_0^y F_1(z)h(z)dz = \frac{1}{l_1}\mathbf{Var}(H(T_{1i})),$$

where

$$H(x) = \int_0^x \frac{[1-F_2(y)]}{f_1(g(y))}dF_1(y) - \int_0^x F_1(g^{-1}(y))dy - \int_0^x F_2(y)dy+$$

$$\int_0^x [1-F_2(y)]dg(y) - \int_0^x [1-F_2(y)]dy = \int_0^x \frac{[1-F_2(y)]}{f_1(g(y))}dF_1(y)-$$

$$F_1(g^{-1}(x))x + \int_0^x y\,dF_1(g^{-1}(y)) + [1-F_2(x)]g(x) + \int_0^x g(y)dF_2(y) - x =$$

$$\int_0^x \frac{[1-F_2(y)]}{f_1(g(y))}dF_1(y) - xF_1(g^{-1}(x)) + g(x)[1-F_2(x)]+$$

Table 1. Significance level of the test

| Sample size | Significance level (%) |
|---|---|
| 50 | 8.47 |
| 100 | 4.63 |
| 170 | 4.37 |
| 200 | 4.43 |
| 400 | 4.77 |

$$\mathbf{E}[\mathbf{1}_{\{g(T_{1i}) \leq x\}} g(T_{1i})] + \mathbf{E}[\mathbf{1}_{\{T_{2i} \leq x\}} g(T_{2i})] - x.$$

□

A consistent estimator of the variance $\sigma^2$ is

$$\hat{\sigma}^2 = \frac{1}{n} \sum_{i=1}^{n} (T_i - \hat{\mu})^2 + \frac{n}{n_1^2} \sum_{i=1}^{n_1} [\hat{H}(T_{1i}) - \hat{\bar{H}}]^2 + \frac{n}{n_2^2} \sum_{i=1}^{n_2} [\hat{Q}(T_{2i}) - \hat{\bar{Q}}]^2,$$

where

$$\hat{H}(x) = \hat{Q}(x) - x\hat{F}_1(\hat{g}^{-1}(x)) + \hat{g}(x)[1 - \hat{F}_2(x)] + \frac{1}{n_1} \sum_{i=1}^{n_1} \mathbf{1}_{\{\hat{g}(T_{1i}) \leq x\}} \hat{g}(T_{1i}) +$$

$$\frac{1}{n_2} \sum_{i=1}^{n_2} \mathbf{1}_{\{T_{2i} \leq x\}} \hat{g}(T_{2i}) - x, \quad \hat{Q}(x) = \frac{1}{n_1} \sum_{i=1}^{n_1} \mathbf{1}_{\{T_{1i} \leq x\}} [1 - \hat{F}_2(T_{1i})]/\hat{f}_1(\hat{g}(T_{1i})),$$

$$\hat{g}^{-1}(x) = \hat{F}_2^{-1}(\hat{F}_1(x)), \quad \hat{\bar{H}} = \frac{1}{n_1} \sum_{i=1}^{n_1} \hat{H}(T_{1i}), \quad \hat{\bar{Q}} = \frac{1}{n_2} \sum_{i=1}^{n_2} \hat{Q}(T_{2i}),$$

the density $f_1$ is estimated by the kernel estimator

$$\hat{f}_1(x) = \frac{1}{n} \sum_{i=1}^{n} \frac{1}{h} K\left(\frac{x - X_{1i}}{h}\right).$$

## 4. Simulations

We did a small simulation study of goodness-of-fit test for the hypothesis $H_0^*$.

Failure times from exponential distribution were simulated:

$$T_{1j} \sim E(\lambda_1), \quad T_{2j} \sim E(\lambda_2), \quad \lambda_2 = r\lambda_1.$$

The distribution function of a redundant system is

$$F(t) = 1 - \frac{\lambda_2 + \lambda_1}{\lambda_1} e^{-\lambda_1 t} + \frac{\lambda_1}{\lambda_2} e^{-(\lambda_1 + \lambda_2)t}.$$

The hypothesis $H_0^*$ is tested using 5 per cent significance level under several sample size n (see Table 1). Number of replications is 3000.

Table 2. Power of the test

| Sample size \ Constant | 0.1 | 0.25 | 0.5 | 0.75 |
|---|---|---|---|---|
| 100 | 6 | 12 | 32 | 85 |
| 170 | 11 | 21 | 79 | 100 |
| 400 | 49 | 88 | 100 | 100 |

Alternative hypothesis $\tilde{H}_0^*$ (at the switching time $y$ the c.d.f. of the stand-by unit has a jump of size $p(1 - F_2(y))$):

$$f_2^{(y)}(x) = f_1\left(x + F_1^{-1}(F_2(y) + p(1 - F_2(y))) - y\right), \quad 0 \leq p \leq 1.$$